\documentclass[12pt,oneside]{amsart}

\usepackage{tikz}
\usepackage{tkz-euclide}

\usepackage{fullpage}
\usepackage[T1]{fontenc}
\usepackage[utf8]{inputenc}
\usepackage{lmodern}
\usepackage{microtype}
\usepackage{amsmath,amssymb,amsthm,mathtools}
\usepackage{booktabs,tabularx,array}
\usepackage{enumitem}
\usepackage{xcolor}
\usepackage[colorlinks=true,linkcolor=blue,citecolor=blue,urlcolor=blue]{hyperref}
\usepackage{todonotes}

\newtheorem{theorem}{Theorem}
\newtheorem{lemma}[theorem]{Lemma}
\newtheorem{proposition}[theorem]{Proposition}
\newtheorem{corollary}[theorem]{Corollary}

\newcommand{\U}{\mathrm U}
\newcommand{\Z}{\mathbb Z}
\newcommand{\N}{\mathbb N}
\newcommand{\R}{\mathbb R}
\newcommand{\CA}{\mathcal A}
\newcommand{\CB}{\mathcal B}
\newcommand{\bits}{\{0,1\}}
\newcommand{\SD}{\operatorname{SD}}
\newcommand{\SIVP}{\operatorname{SIVP}}
\newcommand{\SSP}{\operatorname{SSP}}
\newcommand{\GSIS}{\operatorname{GSIS}}
\newcommand{\NP}{{NP}}
\newcommand{\Set}[2]{\left\{\, #1 \;\middle|\; #2 \,\right\}}

\title{Average-Case Hardness of the Subset Sum Problem
in Groups of Units}
\author{Alexander Ushakov}
\address{Department of Mathematical Sciences,
Stevens Institute of Technology,
Hoboken, NJ 07030}
\email{aushakov@stevens.edu}

\begin{document}

\begin{abstract}
We investigate the worst-case and average-case
complexity of the subset sum problem (SSP)
in finite abelian groups.
We construct a family of composite integers $n_d$
such that solving random planted instances of
SSP in the multiplicative group of units $\U(n_d)$
is at least as hard as approximating the shortest
independent vectors problem on worst-case
$d$-dimensional lattices within a factor of
$O(d\ln^{3/2}d)$.
We further show that SSP over this family of
groups is $\NP$-complete, even when restricted
to instances with a prescribed number of
increments, although positive instances
within this restriction form a strongly
generic set.
Consequently, the decision problem admits a
constant-time heuristic algorithm that is
correct on a strongly generic set of inputs,
whereas the corresponding planted search
problem is hard on average under a
worst-case lattice hardness assumption.
\end{abstract}

\maketitle



\section{Introduction}\label{se:introduction}

The subset sum problem (SSP) in a group $G$
asks, given elements $g_1,\ldots,g_k,g\in G$,
whether there exist
$x_1,\ldots,x_k\in\{0,1\}$ such that
\[
g_1^{x_1}\cdots g_k^{x_k}=g.
\]
The computational complexity of SSP and related
knapsack-type problems in groups has been
investigated extensively; see, for example,
\cite{MNU2015,FNU2016,KLZ2016,LZ2018,
Lohrey2020,NU2020,LRZ2025}.
These works establish polynomial-time
algorithms and hardness results for various
classes of groups, including hyperbolic,
nilpotent, polycyclic, graph, branch,
and finite permutation groups.
They also investigate the behavior of
knapsack-type problems under various
group constructions.

While these works primarily concern worst-case
complexity, the present paper investigates
the average-case complexity of SSP in
finite abelian groups.
Our principal objective is to establish
average-case hardness of SSP in multiplicative
groups of units modulo composite integers,
based on the worst-case hardness of
lattice problems.

More specifically, let $\U(n)$ denote the
multiplicative group of units modulo $n$.
The search version of SSP over $\U(n)$
asks, given elements
$b_1,\ldots,b_k,t\in\U(n)$, to find
$x_1,\ldots,x_k\in\{0,1\}$ satisfying
\[
\prod_{i=1}^k b_i^{x_i}\equiv t\pmod n.
\]
Unlike the classical modular subset sum
problem, which is formulated in the additive
group $\Z/n\Z$, the problem considered here
is formulated in the multiplicative group
$\U(n)$, with group elements represented
by ordinary residue classes modulo $n$.

We study the search problem under a planted
distribution.
The elements $b_1,\ldots,b_k$ are chosen
independently and uniformly from a finite
abelian group $G$, and the target $t$ is
constructed by choosing an independent
uniformly distributed vector
$e\in\{0,1\}^k$ and setting
\[
t=\prod_{i=1}^k b_i^{e_i}.
\]
By construction, every instance admits
at least one solution.
The computational task is to find any
solution, not necessarily the planted
vector $e$.

Our main result establishes average-case
hardness of planted SSP over a suitably
constructed family of multiplicative groups
of units.
More precisely, we construct a family of
squarefree integers $n_d$ and show that
a polynomial-time algorithm solving planted
SSP over $\U(n_d)$ with inverse-polynomial
success probability would yield a randomized
polynomial-time algorithm for approximating
the shortest independent vectors problem
($\SIVP_\gamma$) on worst-case
$d$-dimensional lattices within an
approximation factor
\[
\gamma=O(d\ln^{3/2}d).
\]
Importantly, the hardness result holds even
when the complete prime factorization of
$n_d$ is given as part of the input.

The proof combines a reduction from
Group-SIS to planted SSP with the generalized
worst-case-to-average-case reduction of
Gama, Izabach\`ene, Nguyen, and
Xie~\cite{GINX2016}, building on the
worst-case-to-average-case paradigm
initiated by Ajtai~\cite{Ajtai1996}.

Our first contribution is a reduction
from Group-SIS to planted SSP over an
arbitrary finite abelian group $G$.
Given $k+1$ independent uniformly distributed
elements of $G$, we construct an SSP instance
in which the increments and the target
are independently and uniformly distributed.
We show that the resulting distribution
is statistically close to the planted
distribution, with statistical distance
at most
\[
\frac12\sqrt{\frac{|G|-1}{2^k}}.
\]
Consequently, an algorithm solving
planted SSP with success probability
at least $\delta$ yields a Group-SIS solver
with success probability at least
\[
\delta-\frac12\sqrt{\frac{|G|-1}{2^k}}.
\]
The reduction makes a single call to
the planted SSP solver, performs $O(k)$
additional group operations, and, whenever
successful, produces a nonzero integer
relation of Euclidean norm at most
$\sqrt{k+1}$.

When applied to $\U(n)$, this reduction
requires only modular arithmetic.
In particular, neither the prime factorization
of $n$ nor an explicit decomposition of
$\U(n)$ into cyclic factors is required
to carry out the reduction.

Our second contribution is the construction
of a family of squarefree moduli $n_d$
such that the groups $G_d=\U(n_d)$ satisfy
\[
\operatorname{rank}(G_d)=d,
\qquad
|G_d|=\exp(\Theta(d\ln d)).
\]
The moduli have binary length $O(d\ln d)$,
and both the moduli and their complete
prime factorizations can be computed
in deterministic polynomial time in $d$.
Moreover, the groups $G_d$ admit efficiently
computable cyclic representations.

For a suitable choice of the number
$k_d=\Theta(d\ln d)$ of SSP increments,
the statistical distance between the planted
distribution and the distribution used
in our reduction is bounded above by
$\exp(-\Omega(d\ln d))$.
Thus, the reduction preserves
inverse-polynomial success probabilities.
We then verify that the groups $G_d$
and the chosen parameters satisfy the
hypotheses of the worst-case-to-average-case
reduction of~\cite{GINX2016}, thereby
establishing the claimed connection
with worst-case lattice problems.

Our third contribution concerns the
worst-case and generic-case complexity
of SSP over the same family of groups.
We show that SSP over $G_d$ is
$\NP$-complete, even when the number of
increments is restricted to $k_d$ and
the prime factorization of $n_d$ is given.
At the same time, positive instances
form a strongly generic set with respect
to the stratification by $d$.
Consequently, the constant-time algorithm
that always answers YES is correct on
a strongly generic set of instances.

Thus, the same family of groups exhibits
a striking contrast between the worst-case
complexity of deciding whether a solution
exists, the generic behavior of the
decision problem, and the average-case
complexity of finding a solution.
These results are summarized in
Table~\ref{tab:complexity}.

\begin{table}[ht]
\centering
\begin{tabular}{l|l}
\hline
\textbf{Computational formulation}
  & \textbf{Complexity} \\
\hline
Worst-case decision
  & $\NP$-complete \\
Generic heuristic decision
  & Constant time \\
Average-case planted search
  & Hard under the worst-case \\
  & $\SIVP_\gamma$ assumption \\
\hline
\end{tabular}
\caption{Complexity of $\SSP$ over the family
$G_d=\U(n_d)$ with $k_d$ increments,
where $\gamma=O(d\ln^{3/2}d)$.}
\label{tab:complexity}
\end{table}

Although SSP in multiplicative groups
has natural applications in cryptography,
our primary objective is to understand
its computational complexity.
In particular, the hardness established
here does not rely on the difficulty
of integer factorization or discrete
logarithm computation: both problems
can be solved in deterministic polynomial
time for the constructed family of groups.

The paper is organized as follows.
Section~\ref{se:preliminaries} introduces
the computational problems and probability
distributions used throughout the paper.
Section~\ref{se:reduction} establishes
the reduction from Group-SIS to planted SSP.
Section~\ref{se:groups} constructs the
family of groups of units with the
required structural and computational
properties.
Section~\ref{se:hardness} establishes
average-case hardness of planted SSP
based on the worst-case hardness of
lattice problems.
Finally, Section~\ref{se:ssp} establishes
$\NP$-completeness of SSP over the
constructed family of groups and discusses
its generic-case complexity.

\section{Preliminaries: problems and distributions}
\label{se:preliminaries}

Throughout, $G$ denotes a finite abelian group,
written multiplicatively, with identity element
$1_G$ and order $N=|G|$.
For an integer $n\ge2$, we denote by
$\U(n)=(\Z/n\Z)^\times$ the multiplicative
group of units modulo $n$.

\subsection{Subset sum}

The \emph{Subset Sum Problem} ($\SSP$)
over $G$ is the following decision problem.
Given elements $b_1,\ldots,b_k,t\in G$,
decide whether there exists a vector
$x\in\bits^k$ satisfying
\[
b_1^{x_1}\cdots b_k^{x_k}=t.
\]

\subsection{Planted subset sum}

Fix an integer $k\ge1$.
The \emph{planted subset sum distribution}
$\mathcal P_{G,k}$ on $G^{k+1}$ is defined
by the following random experiment:
\begin{enumerate}
\item Sample $b_1,\ldots,b_k$ independently
and uniformly from $G$.
\item Sample a vector $e=(e_1,\ldots,e_k)$
uniformly from $\bits^k$, independently
of $b_1,\ldots,b_k$.
\item Set $t=b_1^{e_1}\cdots b_k^{e_k}$.
\item Output $(b_1,\ldots,b_k,t)$.
\end{enumerate}

The \emph{Planted Subset Sum Problem}
$\SSP_{G,k}$ is the search problem of finding,
given an instance
$(b_1,\ldots,b_k,t)$ sampled from
$\mathcal P_{G,k}$,
a vector $x\in\bits^k$ satisfying
\[
b_1^{x_1}\cdots b_k^{x_k}=t.
\]
By construction, every instance sampled
from $\mathcal P_{G,k}$ admits at least
one solution, namely the planted vector $e$.
The objective, however, is to find any
solution, not necessarily $e$.

\subsection{Probabilistic and asymptotic notation}

For a finite set $S$, the notation
$x\leftarrow S$ means that $x$ is sampled
uniformly from $S$.
Multiple samples are assumed to be independent
unless otherwise specified.
For a randomized algorithm $\CA$,
the notation $x\leftarrow\CA(y)$ denotes
the output of $\CA$ on input $y$.

We write $f(d)\ge1/\operatorname{poly}(d)$
if there exists a polynomial $p(d)$,
positive for all sufficiently large $d$,
such that $f(d)\ge1/p(d)$ for all
sufficiently large $d$.

Throughout, randomized polynomial-time
algorithms are required to run in polynomial
time on every input, including inputs outside
the support of the distribution under
consideration.

\subsection{Exact Cover by 3-Sets}

The \emph{Exact Cover by 3-Sets problem}
(X3C) is the following decision problem.
Given a finite set $S=\{s_1,\ldots,s_m\}$,
where $m$ is divisible by $3$, and a collection
$S_1,\ldots,S_\ell$ of three-element subsets
of $S$, decide whether there exists a
subcollection that partitions $S$.
Equivalently, the problem is to decide
whether there exists a vector $x\in\bits^\ell$
such that
\[
|\Set{i}{1\le i\le\ell,\ x_i=1,\ s_j\in S_i}|=1
\qquad \text{for every } j=1,\ldots,m.
\]
It is well known that X3C is $\NP$-complete
(see~\cite{GareyJohnson1979}).

\subsection{Shortest Independent Vectors Problem}

Let $L\subseteq\R^d$ be a full-rank lattice.
Its $d$th \emph{successive minimum}
$\lambda_d(L)$ is the smallest real number
$r$ such that $L$ contains $d$ linearly
independent vectors of Euclidean norm
at most $r$.

For an \emph{approximation factor} $\gamma\ge1$,
the \emph{Shortest Independent Vectors Problem}
$\SIVP_\gamma$ is the worst-case search
problem of finding, given a basis of $L$,
linearly independent vectors
$v_1,\ldots,v_d\in L$ satisfying
\[
\max_{1\le i\le d}\|v_i\|_2
\le\gamma\lambda_d(L).
\]
For $\epsilon>0$, the
\emph{smoothing parameter} $\eta_\epsilon(L)$
is defined by
\[
\eta_\epsilon(L)
=
\inf\Set{s>0}{
\sum_{y\in L^*\setminus\{0\}}
e^{-\pi s^2\|y\|_2^2}
\le\epsilon
},
\]
where $L^*$ denotes the dual lattice.
Thus, $\eta_\epsilon(L)$ is the smallest
scale at which the total Gaussian mass
of nonzero dual lattice vectors is at most
$\epsilon$. In particular, decreasing
$\epsilon$ increases the smoothing parameter.
By~\cite[Lemma~3.3]{MR2007}, for $0<\epsilon<1$,
we have
\[
\eta_\epsilon(\Z^d)
\le
\sqrt{\frac{\ln(2d(1+1/\epsilon))}{\pi}}.
\]
Consequently, for
$\epsilon(d)=\exp(-(\ln d)^2)$,
\[
\eta_{\epsilon(d)}(\Z^d)=O(\ln d).
\]

\subsection{The Group-SIS problem}

Let $M\ge1$ and $\beta>0$.
The \emph{Group Short Integer Solution problem}
(Group-SIS), introduced by Micciancio
under the name \emph{Homogeneous SIS}
in~\cite{Micciancio2004} and studied
in~\cite{GINX2016}, is the randomized
search problem $\GSIS_{G,M,\beta}$ of finding,
given independent uniformly distributed
elements $g_1,\ldots,g_M\in G$,
a nonzero vector $z\in\Z^M$ satisfying
\[
g_1^{z_1}\cdots g_M^{z_M}=1_G,
\qquad
\|z\|_2\le\beta.
\]
As noted in~\cite[Section~3.2]{GINX2016},
a solution is guaranteed to exist whenever
$\beta\ge\sqrt{M}N^{1/M}$.
This follows from the pigeonhole principle
applied to integer vectors of bounded norm.

\subsection{Worst-case hardness of Group-SIS}

We recall the worst-case-to-average-case
reduction of Gama, Izabach\`ene, Nguyen,
and Xie~\cite{GINX2016}.

A finite abelian group $G$ is called
\emph{explicit} if a decomposition
\[
G\cong C_{q_1}\times\cdots\times C_{q_r},
\]
where $C_m$ denotes the cyclic group of order $m$,
is known together with a polynomial-time
computable isomorphism from the product
of cyclic groups to $G$.
Here $r=\operatorname{rank}(G)$ denotes
the minimum number of generators of $G$.
If the inverse isomorphism is also
computable in polynomial time,
the group is called \emph{fully explicit}.

A function $\epsilon:\N\to\R_{>0}$
is called \emph{negligible} if it tends
to zero faster than any inverse polynomial,
that is, if
\[
\forall c\in\N_{>0}\ 
\exists d_0\in\N\
\forall d\ge d_0:\quad
\epsilon(d)<d^{-c}.
\]
The following theorem is
\cite[Corollary~5.2]{GINX2016},
restated in our notation.

\begin{theorem}\label{thm:GINX}
Let $\epsilon:\N\to\R_{>0}$ be negligible,
and let $(H_d)$ be a family of explicit
finite abelian groups with
$r_d=\operatorname{rank}(H_d)$.
Suppose that $M_d\ge1$ and $\beta_d>0$
satisfy
\[
d^{r_d}
\left(
\eta_{\epsilon(d)}(\Z^d)
\sqrt{\frac{2d}{\pi}}\,\beta_d
\right)^{\max(d,r_d)}
\le
|H_d|
\le
\left(\frac{\beta_d}{\sqrt{M_d}}\right)^{M_d}.
\]
Suppose that there exists a randomized
polynomial-time algorithm $\CA$ such that
\[
\Pr_{\substack{
g_1,\ldots,g_{M_d}\leftarrow H_d\\
z\leftarrow\CA(g_1,\ldots,g_{M_d})
}}
\left[
z\in\Z^{M_d}\setminus\{0\},\
\|z\|_2\le\beta_d,\
\prod_{i=1}^{M_d}g_i^{z_i}=1_{H_d}
\right]
\ge\frac{1}{\operatorname{poly}(d)}.
\]
Then, using polynomially many calls
to $\CA$, one can solve worst-case
$d$-dimensional $\SIVP_{\gamma_d}$
in randomized polynomial time, where
\[
\gamma_d
=
\eta_{\epsilon(d)}(\Z^d)
\sqrt{\frac{d}{\pi}}\,\beta_d.
\]
\end{theorem}

We will apply Theorem~\ref{thm:GINX}
to a family of groups of units
$G_d=\U(n_d)$ constructed below.

\section{Reducing Group-SIS to planted subset sum}
\label{se:reduction}

We establish a reduction from Group-SIS
to planted SSP over an arbitrary finite
abelian group $G$.
The reduction relies on the statistical
closeness of planted and independent-target
distributions and does not require any
lattice hardness assumption.

\subsection{Subset-sum hashing}

For $b=(b_1,\ldots,b_k)\in G^k$, define
\[
h_b:\bits^k\longrightarrow G,
\qquad
h_b(x)=b_1^{x_1}\cdots b_k^{x_k}.
\]

The following lemma shows that these
functions form a universal family of
hash functions.

\begin{lemma}\label{lem:universal}
For any distinct $x,y\in\bits^k$,
\[
\Pr_{b\leftarrow G^k}
[h_b(x)=h_b(y)]
=
\frac1N.
\]
\end{lemma}

\begin{proof}
Choose an index $j$ such that $x_j\ne y_j$.
The equality $h_b(x)=h_b(y)$ is equivalent to
\[
b_j^{x_j-y_j}
=
\prod_{i\ne j}b_i^{y_i-x_i}.
\]
Conditioning on the elements $b_i$ with
$i\ne j$ fixes the right-hand side.
Since $x_j-y_j\in\{-1,1\}$,
exactly one value of $b_j$ satisfies
this equation.
As $b_j$ is uniformly distributed
and independent of the other elements,
the probability is $1/N$.
\end{proof}

\subsection{Planted and independent targets}

For distributions $P,Q$ on a finite set
$\Omega$, their \emph{statistical distance}
is defined by
\[
\SD(P,Q)
=
\frac12\sum_{\omega\in\Omega}
|P(\omega)-Q(\omega)|.
\]
Let $\mathcal U_{G,k}$ denote the
\emph{independent-target distribution}
on $G^{k+1}$, obtained by sampling
$b_1,\ldots,b_k,t$ independently
and uniformly from $G$.
The following estimate is an application
of the leftover hash lemma.

\begin{lemma}\label{lem:lhl}
For every $k\ge1$,
\[
\SD(\mathcal P_{G,k},\mathcal U_{G,k})
\le
\frac12\sqrt{\frac{N-1}{2^k}}.
\]
\end{lemma}

\begin{proof}
For fixed $b\in G^k$, let $P_b$ denote
the distribution of $h_b(X)$ for
uniformly distributed $X\in\bits^k$,
and let $U_G$ denote the uniform
distribution on $G$.
Under both $\mathcal P_{G,k}$ and
$\mathcal U_{G,k}$, the vector $b$
is uniformly distributed over $G^k$.
Hence,
\[
\SD(\mathcal P_{G,k},\mathcal U_{G,k})
=
\mathbb E_b[\SD(P_b,U_G)],
\]
where the expectation is taken over
uniformly distributed $b\in G^k$.

By Cauchy--Schwarz and Jensen's inequality,
\[
\SD(\mathcal P_{G,k},\mathcal U_{G,k})
\le
\frac12\mathbb E_b
\sqrt{
N\sum_{t\in G}
\left(P_b(t)-\frac1N\right)^2
}\le
\frac12\sqrt{
N\mathbb E_b
\sum_{t\in G}
\left(P_b(t)-\frac1N\right)^2
}.
\]

Let $X,Y$ be independent uniform
elements of $\bits^k$.
Then
\begin{align*}
\sum_{t\in G}P_b(t)^2
&=
\sum_{t\in G}
\Pr[h_b(X)=t]\Pr[h_b(Y)=t]\\
&=
\sum_{t\in G}
\Pr[h_b(X)=t,\ h_b(Y)=t]\\
&=
\Pr[h_b(X)=h_b(Y)].
\end{align*}
By Lemma~\ref{lem:universal},
\[
\mathbb E_b\sum_{t\in G}P_b(t)^2
=
\Pr[X=Y]+\frac1N\Pr[X\ne Y]
=
2^{-k}+\frac{1-2^{-k}}N.
\]
Consequently,
\[
\mathbb E_b\sum_{t\in G}
\left(P_b(t)-\frac1N\right)^2
=
\mathbb E_b\sum_{t\in G}P_b(t)^2
-\frac1N
=
\frac{N-1}{N2^k}.
\]
Substituting this estimate into
the preceding inequality proves
the lemma.
\end{proof}

\begin{corollary}\label{cor:solvability}
An instance sampled from
$\mathcal U_{G,k}$ has a solution
with probability at least
\[
1-\frac12\sqrt{\frac{N-1}{2^k}}.
\]
\end{corollary}

\begin{proof}
Every instance sampled from
$\mathcal P_{G,k}$ has a solution.
The result follows from
Lemma~\ref{lem:lhl}, since statistical
distance bounds the difference in
probabilities of any event.
\end{proof}

Consequently, for any family of groups
and parameters for which $N/2^k$
is negligible in the security parameter,
independent-target SSP instances are
solvable with overwhelming probability.

\subsection{The reduction}

We now show how to use an algorithm for
planted SSP to solve Group-SIS.

\begin{proposition}\label{prop:reduction}
Let $G$ be a finite abelian group of order $N$,
and let $k\ge1$.
Suppose that a randomized algorithm $\CA$
solves $\SSP_{G,k}$ with probability
at least $\delta$, taken over the planted
distribution and the internal randomness
of $\CA$.

Then there exists a randomized algorithm
$\CB$ for $\GSIS_{G,k+1,\sqrt{k+1}}$
that outputs a valid solution with probability
at least
\[
\delta-\frac12\sqrt{\frac{N-1}{2^k}},
\]
and otherwise reports failure.
The algorithm $\CB$ makes a single call
to $\CA$ and performs $O(k)$ additional
group operations.
\end{proposition}

\begin{proof}
Given independent uniformly distributed
elements $g_1,\ldots,g_{k+1}\in G$,
the algorithm $\CB$ proceeds as follows:
\begin{enumerate}
\item
Set $b_i=g_i$ for $1\le i\le k$
and $t=g_{k+1}^{-1}$.
\item
Run $\CA$ on the instance
$(b_1,\ldots,b_k,t)$.
\item
If $\CA$ returns a vector $x\in\bits^k$
satisfying $h_b(x)=t$, output
$z=(x_1,\ldots,x_k,1)$.
Otherwise, report failure.
\end{enumerate}

We first verify the correctness of $\CB$.
Whenever $\CB$ outputs a vector $z$, we have
\[
\prod_{i=1}^{k+1}g_i^{z_i}
=
h_b(x)g_{k+1}
=
tg_{k+1}
=
1_G.
\]
Moreover, $z\ne0$ and
\[
\|z\|_2^2
=
1+\sum_{i=1}^k x_i
\le k+1.
\]
Thus every vector returned by $\CB$
is a valid solution of
$\GSIS_{G,k+1,\sqrt{k+1}}$.

It remains to estimate the success probability
of $\CB$.
Since inversion preserves the uniform
distribution on $G$, the instance
$(b_1,\ldots,b_k,t)$ supplied to $\CA$
has distribution $\mathcal U_{G,k}$.
By Lemma~\ref{lem:lhl},
\[
\SD(\mathcal P_{G,k},\mathcal U_{G,k})
\le\frac12\sqrt{\frac{N-1}{2^k}}.
\]
Consequently, the success probability of
$\CA$ under $\mathcal U_{G,k}$ differs
from its success probability under
$\mathcal P_{G,k}$ by at most this quantity.
Since $\CB$ succeeds precisely when
$\CA$ returns a valid SSP solution,
the success probability of $\CB$ is at least
\[
\delta-\frac12\sqrt{\frac{N-1}{2^k}}.
\]

Finally, apart from a single call to $\CA$,
the algorithm $\CB$ performs one inversion
and at most $k$ multiplications in $G$
to construct the instance and verify
the returned solution.
\end{proof}

For $G=\U(n)$, the reduction in
Proposition~\ref{prop:reduction} can be
implemented entirely using modular
arithmetic. In particular, it requires
neither the factorization of $n$ nor
an explicit decomposition of $\U(n)$
into cyclic factors.
The factorization and primitive roots
introduced in Section~\ref{se:groups}
are needed only to establish the structural
properties of the groups $G_d$ required
for the worst-case-to-average-case reduction
of~\cite{GINX2016}.

\begin{corollary}\label{cor:density}
Suppose that a randomized algorithm $\CA$
solves $\SSP_{G,k}$ with probability
at least $\delta>0$.
If
\[
k\ge\log_2 N+2\log_2(1/\delta),
\]
then there exists a randomized algorithm
$\CB$ for $\GSIS_{G,k+1,\sqrt{k+1}}$
with success probability at least $\delta/2$.
The algorithm $\CB$ makes a single call
to $\CA$ and performs $O(k)$ additional
group operations.
\end{corollary}

\begin{proof}
The assumption on $k$ implies that
$2^k\ge N/\delta^2$, and hence
\[
\frac12\sqrt{\frac{N-1}{2^k}}
\le\frac12\sqrt{\frac{N}{2^k}}
\le\frac{\delta}{2}.
\]
The result now follows directly from
Proposition~\ref{prop:reduction}.
\end{proof}

\section{Groups of units with controlled parameters}
\label{se:groups}

We specialize the preceding reduction to
multiplicative groups of units modulo composite
integers. We construct a family of groups
$G_d=\U(n_d)$ of rank $d$ and order
$\exp(\Theta(d\ln d))$, whose elements
can be represented using $O(d\ln d)$ bits.

Each modulus $n_d$ is a product of $d$
distinct primes of comparable size.
The Chinese remainder theorem gives
a decomposition of $G_d$ into cyclic factors.
The parameters are chosen so that the
statistical distance in Lemma~\ref{lem:lhl}
is negligible.
Both the moduli $n_d$ and their prime
factorizations can be computed in
deterministic polynomial time from $d$.

\subsection{Moduli $n_d$}\label{se:moduli}

We begin by establishing the existence of
sufficiently many primes in the interval
$[d^4,2d^4]$.

\begin{lemma}\label{le:primes}
For every sufficiently large integer $d$,
the interval $[d^4,2d^4]$ contains at least
$d$ distinct odd primes.
\end{lemma}

\begin{proof}
Let $\pi(x)$ denote the number of primes
not exceeding $x$.
By the prime number theorem,
\[
\pi(2d^4)-\pi(d^4)
\sim
\frac{d^4}{4\ln d}
\qquad\text{as }d\to\infty.
\]
Since
\[
\lim_{d\to\infty}
\frac{d^4}{4d\ln d}
=\infty,
\]
we have $\pi(2d^4)-\pi(d^4)\geq d$
for all sufficiently large $d$.
Moreover, all primes in the interval
$[d^4,2d^4]$ are odd for $d\geq2$.
This proves the claim.
\end{proof}

For each sufficiently large $d$, let
$p_{d,1}<\cdots<p_{d,d}$ be the first
$d$ primes in $[d^4,2d^4]$, and define
\begin{equation}\label{eq:n_d}
n_d=\prod_{j=1}^d p_{d,j},
\qquad
G_d=\U(n_d).
\end{equation}

\begin{lemma}\label{le:moduli-computation}
The modulus $n_d$, together with its complete
prime factorization, can be computed in
deterministic polynomial time in $d$.
\end{lemma}

\begin{proof}
Enumerate the integers in the interval
$[d^4,2d^4]$ in increasing order and test
each for primality until $d$ primes have
been found.
By Lemma~\ref{le:primes}, this procedure
terminates for all sufficiently large $d$.

The interval contains $O(d^4)$ integers,
each of bit length $O(\log d)$.
Since every candidate is bounded by $2d^4$,
its primality can be tested by trial division
using $O(d^2)$ arithmetic operations on
integers of bit length $O(\log d)$.
Thus, the required primes can be found in
deterministic polynomial time in $d$.

Finally, multiplying the resulting $d$ primes
yields the modulus $n_d$, whose bit length
is $O(d\log d)$.
Hence both $n_d$ and its complete prime
factorization can be computed in
deterministic polynomial time in $d$.
\end{proof}

We next determine the structure and size
of the groups $G_d$.

\begin{proposition}\label{prop:group-parameters}
For every sufficiently large $d$, the group
$G_d=\U(n_d)$ satisfies:
\begin{itemize}
\item $\operatorname{rank}(G_d)=d$;
\item $|G_d|=\exp(\Theta(d\ln d))$;
\item $\ln n_d=\Theta(d\ln d)$.
\end{itemize}
In particular, elements of $G_d$ can be
represented using $O(d\ln d)$ bits.
\end{proposition}

\begin{proof}
Since $n_d$ is squarefree, the Chinese
remainder theorem gives
\[
G_d
\cong
\prod_{j=1}^d\U(p_{d,j})
\cong
\prod_{j=1}^d C_{p_{d,j}-1},
\]
where $C_m$ denotes the cyclic group
of order $m$.
Each $p_{d,j}$ is odd, so
$p_{d,j}-1$ is even. Consequently,
\[
G_d/G_d^2\cong(\Z/2\Z)^d,
\]
where $G_d^2$ denotes the subgroup
of squares in $G_d$.
Thus $G_d$ requires at least $d$
generators, while the displayed
cyclic decomposition shows that
$d$ generators suffice.
Hence $\operatorname{rank}(G_d)=d$.

Furthermore,
\[
|G_d|=\prod_{j=1}^d(p_{d,j}-1),
\]
and therefore
\[
(d^4-1)^d
\le |G_d|
\le (2d^4)^d.
\]
Similarly,
\[
d^{4d}\le n_d\le(2d^4)^d.
\]
Taking logarithms yields
\[
\ln|G_d|=\Theta(d\ln d),
\qquad
\ln n_d=\Theta(d\ln d).
\]
Finally, each element of $G_d$ can be
represented by an integer in
$\{0,\ldots,n_d-1\}$, requiring
$O(\ln n_d)=O(d\ln d)$ bits.
\end{proof}

\subsection{Reduction over groups of units}

For every sufficiently large $d$, define
\[
M_d=\lceil12d\ln d\rceil,
\qquad
k_d=M_d-1,
\qquad
\varepsilon_d=
\frac12\sqrt{\frac{|G_d|-1}{2^{k_d}}}.
\]
The following lemma shows that the statistical
distance between the planted and independent-target
distributions is negligible for these parameters.

\begin{lemma}\label{lem:units-distance}
$\varepsilon_d=\exp(-\Omega(d\ln d))$.
\end{lemma}

\begin{proof}
By Proposition~\ref{prop:group-parameters},
\[
\ln|G_d|\le4d\ln d+d\ln2.
\]
On the other hand,
$k_d\ge12d\ln d-1$.
Consequently,
\[
\begin{aligned}
\ln\left(\frac{|G_d|}{2^{k_d}}\right)
&\le
4d\ln d+d\ln2
-(12d\ln d-1)\ln2\\
&=
-(12\ln2-4)d\ln d
+(d+1)\ln2.
\end{aligned}
\]
Since $12\ln2>4$, the right-hand side
is $-\Omega(d\ln d)$.
Therefore,
\[
\varepsilon_d
\le
\frac12\sqrt{\frac{|G_d|}{2^{k_d}}}
=
\exp(-\Omega(d\ln d)).
\]
\end{proof}

\begin{theorem}\label{thm:units}
Suppose that a randomized algorithm $\CA$
solves $\SSP_{G_d,k_d}$ with probability
at least $\delta_d$.
Then there exists a randomized algorithm
$\CB$ for $\GSIS_{G_d,M_d,\sqrt{M_d}}$
with success probability at least
\[
\delta_d-\varepsilon_d,
\qquad
\text{where }
\varepsilon_d=\exp(-\Omega(d\ln d)).
\]
The algorithm $\CB$ makes a single call
to $\CA$ and has additional running time
polynomial in $d$.
\end{theorem}

\begin{proof}
Apply Proposition~\ref{prop:reduction}
with $G=G_d$ and $k=k_d$.
Since $M_d=k_d+1$, the resulting algorithm
$\CB$ solves $\GSIS_{G_d,M_d,\sqrt{M_d}}$
with probability at least
$\delta_d-\varepsilon_d$, where
\[
\varepsilon_d=
\frac12\sqrt{\frac{|G_d|-1}{2^{k_d}}}.
\]
By Lemma~\ref{lem:units-distance},
we have $\varepsilon_d=\exp(-\Omega(d\ln d))$.

Apart from a single call to $\CA$,
the algorithm $\CB$ performs $O(k_d)$
group operations in $G_d=\U(n_d)$.
Since $k_d=O(d\ln d)$ and
$\ln n_d=O(d\ln d)$, and since
multiplication and inversion modulo $n_d$
can be performed in polynomial time in
$\ln n_d$, the computational overhead
of the reduction is polynomial in $d$.
\end{proof}

\begin{table}[ht]
\centering
\small
\begin{tabular}{c|l|r|r}
\hline
$d$ & Primes $p_{d,1},\ldots,p_{d,d}$
    & $n_d$ & $k_d$ \\
\hline
2 & \scalebox{0.8}{$17,19$}
  & \scalebox{0.8}{$323$} & 16 \\
3 & \scalebox{0.8}{$83,89,97$}
  & \scalebox{0.8}{$716539$} & 39 \\
4 & \scalebox{0.8}{$257,263,269,271$}
  & \scalebox{0.8}{$4927316309$} & 66 \\
5 & \scalebox{0.8}{$631,641,643,647,653$}
  & \scalebox{0.8}{$109879284718823$} & 96 \\
6 & \scalebox{0.8}{$1297,1301,1303,1307,1319,1321$}
  & \scalebox{0.8}{$5007084136217062463$} & 129 \\
7 & \scalebox{0.8}{$2411,2417,2423,2437,2441,2447,2459$}
  & \scalebox{0.8}{$505409128389380177323741$} & 163 \\
8 & \scalebox{0.8}{$4099,4111,4127,4129,4133,4139,4153,4157$}
  & \scalebox{0.8}{$84802281027518993007807702649$} & 199 \\
\hline
\end{tabular}
\caption{The moduli $n_d$, constructed using
the first $d$ primes in $[d^4,2d^4]$,
and the corresponding subset-sum
parameter $k_d=\lceil12d\ln d\rceil-1$.}
\label{tab:moduli}
\end{table}

\section{Average-case hardness of planted subset sum}
\label{se:hardness}

We now combine Theorem~\ref{thm:units}
with the worst-case-to-average-case
reduction of Gama, Izabach\`ene, Nguyen,
and Xie~\cite{GINX2016}.
The main result of this section shows
that solving random planted instances
of SSP in the groups $G_d=\U(n_d)$
is at least as hard as approximating
SIVP on worst-case $d$-dimensional
lattices within a factor of
$O(d\ln^{3/2}d)$.
We first verify that these groups admit
the explicit representations required
by~\cite{GINX2016}.

\subsection{Explicit cyclic representations}

Although the cyclic decomposition of $G_d$
follows immediately from the Chinese remainder
theorem, we need to establish that the
corresponding isomorphism and its inverse
can be computed efficiently.
The following lemma shows that this is
indeed the case for our family of groups.

\begin{lemma}\label{lem:explicit-units}
Let $G_d=\U(n_d)$ be as in~\eqref{eq:n_d}.
One can compute primitive roots $a_{d,j}$
modulo $p_{d,j}$ in deterministic polynomial
time in $d$.
Moreover, the resulting isomorphism
\[
\Phi_d:
\prod_{j=1}^d
\Z/(p_{d,j}-1)\Z
\longrightarrow G_d
\]
and its inverse can both be evaluated
in deterministic polynomial time in $d$.
\end{lemma}

\begin{proof}
By Lemma~\ref{le:moduli-computation},
the prime factors $p_{d,j}$ of $n_d$
can be computed in deterministic polynomial
time in $d$.
Since $p_{d,j}\le2d^4$, the integer
$p_{d,j}-1$ can be factored by trial
division in polynomial time in $d$.
Using this factorization, one can test
whether an element of $\U(p_{d,j})$
is a primitive root by modular
exponentiation.
Since a primitive root exists, exhaustive
search finds one after at most
$p_{d,j}-1=O(d^4)$ tests.

The isomorphism $\Phi_d$ maps a tuple
$(x_1,\ldots,x_d)$ to the unique element
$u\in G_d$ satisfying
\[
u\equiv a_{d,j}^{x_j}\pmod{p_{d,j}},
\qquad j=1,\ldots,d.
\]
It can therefore be evaluated in polynomial
time using modular exponentiation and
the Chinese remainder theorem.

To compute its inverse, reduce an element
$u\in G_d$ modulo each $p_{d,j}$
and find $0\le x_j<p_{d,j}-1$ satisfying
\[
a_{d,j}^{x_j}\equiv u\pmod{p_{d,j}}.
\]
Since $p_{d,j}=O(d^4)$, each discrete
logarithm can be computed by exhaustive
search in polynomial time in $d$.
Thus both $\Phi_d$ and $\Phi_d^{-1}$
are computable in deterministic polynomial
time in $d$.
\end{proof}

The polynomial-time computability of
$\Phi_d^{-1}$ relies on the fact that
$p_{d,j}=O(d^4)$.
For general groups of units, computing
cyclic coordinates may require solving
discrete logarithm problems of
substantially greater complexity.
The explicit cyclic representation
is needed to apply the
worst-case-to-average-case reduction
of~\cite{GINX2016}.
By contrast, the reduction of
Proposition~\ref{prop:reduction}
operates directly on residue classes
and does not require this representation.

\subsection{The main hardness theorem}

Having established that the groups $G_d$
admit explicit cyclic representations,
we now show that our choice of parameters
satisfies the numerical conditions of
Theorem~\ref{thm:GINX}.
Combining this theorem with
Theorem~\ref{thm:units}, we obtain a
reduction from worst-case SIVP to
planted SSP.

\begin{theorem}\label{thm:hardness}
Let $(G_d)$ be the family of groups
$G_d=\U(n_d)$ defined in~\eqref{eq:n_d}.
Set
\[
M_d=\lceil12d\ln d\rceil,
\qquad
k_d=M_d-1,
\qquad
\beta_d=2\sqrt{M_d},
\]
and let
\[
\epsilon_d=\exp(-(\ln d)^2).
\]
Suppose that there exists a randomized
polynomial-time algorithm $\CA$ such that
\[
\Pr_{\substack{
b_1,\ldots,b_{k_d}\leftarrow G_d\\
e\leftarrow\bits^{k_d}\\
x\leftarrow\CA(n_d,b_1,\ldots,b_{k_d},h_b(e))
}}
\left[
x\in\bits^{k_d},\
h_b(x)=h_b(e)
\right]
\ge\frac{1}{\operatorname{poly}(d)}.
\]
Then, using polynomially many calls
to $\CA$, one can solve worst-case
$d$-dimensional $\SIVP_{\gamma_d}$
in randomized polynomial time, where
\[
\gamma_d
=
\eta_{\epsilon_d}(\Z^d)
\sqrt{\frac{d}{\pi}}\,\beta_d
=
O(d\ln^{3/2}d).
\]
\end{theorem}

\begin{proof}
We apply Theorem~\ref{thm:GINX}
with $H_d=G_d$, $r_d=d$, and
$\epsilon(d)=\epsilon_d$.
Note that
\[
\epsilon_d=\exp(-(\ln d)^2)=d^{-\ln d}
\]
is negligible.
By Lemma~\ref{lem:explicit-units} and
Proposition~\ref{prop:group-parameters},
the group $G_d$ has rank $d$ and admits
an explicit cyclic representation
computable in polynomial time.
We verify the two numerical conditions
of Theorem~\ref{thm:GINX}.

The first condition is
\[
|G_d|
\le
\left(
\frac{\beta_d}{\sqrt{M_d}}
\right)^{M_d}
=
2^{M_d}.
\]
By Proposition~\ref{prop:group-parameters},
\[
\ln|G_d|
\le4d\ln d+d\ln2,
\]
whereas
\[
M_d\ln2
=
12(\ln2)d\ln d+O(1).
\]
Since $12\ln2>4$, the first
condition holds for all sufficiently
large $d$.

The second condition is
\[
|G_d|
\ge
d^d
\left(
\eta_{\epsilon_d}(\Z^d)
\sqrt{\frac{2d}{\pi}}\,
\beta_d
\right)^d.
\]
Recall the standard bound
(see~\cite[Lemma~3.3]{MR2007})
\[
\eta_\epsilon(\Z^d)
\le
\sqrt{
\frac{\ln(2d(1+1/\epsilon))}{\pi}
}.
\]
Since
\[
\epsilon_d=\exp(-(\ln d)^2),
\]
we obtain
\[
\eta_{\epsilon_d}(\Z^d)=O(\ln d).
\]
Together with
\[
\beta_d=O(\sqrt{d\ln d}),
\]
this gives
\[
d\,
\eta_{\epsilon_d}(\Z^d)
\sqrt{\frac{2d}{\pi}}\,
\beta_d
=
O(d^2\ln^{3/2}d).
\]

On the other hand, by
Proposition~\ref{prop:group-parameters},
\[
|G_d|\ge(d^4-1)^d.
\]
Since
\[
d^4-1
\ge
C d^2\ln^{3/2}d
\]
for every fixed $C>0$ and all
sufficiently large $d$,
the second condition also holds.
Thus, both numerical conditions of
Theorem~\ref{thm:GINX} are satisfied
for all sufficiently large $d$.

Let $\delta_d$ denote the success
probability of $\CA$ over the planted
distribution and its internal randomness.
By assumption,
\[
\delta_d\ge\frac{1}{\operatorname{poly}(d)}.
\]
By Theorem~\ref{thm:units}, we obtain
a randomized polynomial-time algorithm
for $\GSIS_{G_d,M_d,\sqrt{M_d}}$
with success probability at least
\[
\delta_d-\varepsilon_d
\ge\frac{\delta_d}{2}
\]
for all sufficiently large $d$, since
$\delta_d$ is bounded below by an inverse
polynomial and
$\varepsilon_d=\exp(-\Omega(d\ln d))$.
Since
\[
\sqrt{M_d}\le\beta_d,
\]
the same algorithm solves
$\GSIS_{G_d,M_d,\beta_d}$
with inverse-polynomial success
probability.

Theorem~\ref{thm:GINX} therefore yields
a randomized polynomial-time algorithm
for worst-case $d$-dimensional
$\SIVP_{\gamma_d}$, where
\[
\gamma_d
=
\eta_{\epsilon_d}(\Z^d)
\sqrt{\frac{d}{\pi}}\,\beta_d
=
O(d\ln^{3/2}d).
\]
This proves the theorem.
\end{proof}

By Lemma~\ref{le:moduli-computation},
the complete prime factorization of $n_d$
can be computed in deterministic polynomial
time in $d$.
Thus, providing the prime factorization
of $n_d$ as part of the input does not
change the computational complexity
of the planted subset sum problem
over the family $(G_d)$.
Moreover, by Lemma~\ref{lem:explicit-units},
discrete logarithms in the cyclic factors
of $G_d$ can also be computed in
deterministic polynomial time in $d$.
Hence, the conditional average-case hardness
established in Theorem~\ref{thm:hardness}
relies on neither the hardness of integer
factorization nor that of discrete
logarithm computation.

\section{$\SSP$ over $\U(n_d)$ is $\NP$-complete}
\label{se:ssp}

In this section, we establish the worst-case
hardness of the decision version of the
subset sum problem over the family of groups
$G_d=\U(n_d)$ constructed above.
In fact, the problem remains $\NP$-complete
even when the number of increments is
restricted to $k_d=\lceil12d\ln d\rceil-1$
and the prime factorization of $n_d$
is included in the input.

\begin{theorem}\label{thm:ssp-np-complete}
The problem of deciding, given the modulus $n_d$
and integers $b_1,\ldots,b_{k_d},t$ coprime
to $n_d$, whether there exists a vector
$x\in\bits^{k_d}$ satisfying
\[
b_1^{x_1}\cdots b_{k_d}^{x_{k_d}}
\equiv t\pmod{n_d}
\]
is $\NP$-complete.
\end{theorem}

\begin{proof}
Membership in $\NP$ is immediate, since
a certificate $x\in\bits^{k_d}$ can be
verified in polynomial time using
modular arithmetic.

To establish $\NP$-hardness, we construct
a polynomial-time reduction from
the Exact Cover by 3-Sets problem (X3C).
Let $S=\{s_1,\ldots,s_m\}$ be a finite set,
and let $S_1,\ldots,S_\ell$ be a collection
of three-element subsets of $S$.
The problem is to decide whether there
exists a subcollection of these subsets
that partitions $S$.

Fix a constant $d_0$ sufficiently large
for Lemma~\ref{le:primes} to hold, and set
$d=\max(d_0,m,\ell)$.
In particular, we have $k_d\ge\ell$.
Recall that
\[
n_d=\prod_{j=1}^d p_{d,j},
\qquad
d^4\le p_{d,j}\le2d^4,
\]
where the $p_{d,j}$ are distinct primes.

For each $1\le j\le m$, let $a_j$ be
a primitive root modulo $p_{d,j}$.
Since $p_{d,j}=O(d^4)$, such elements
can be computed by exhaustive search
in deterministic polynomial time in $d$.

For each $1\le i\le\ell$, define
an integer $b_i$ by the congruences
\[
b_i\equiv
\begin{cases}
a_j \pmod{p_{d,j}},
    & 1\le j\le m,\ s_j\in S_i,\\
1 \pmod{p_{d,j}},
    & \text{otherwise}.
\end{cases}
\]
Set $b_i=1$ for $\ell<i\le k_d$.
Similarly, define the target $t$ by
\[
t\equiv
\begin{cases}
a_j \pmod{p_{d,j}},
    & 1\le j\le m,\\
1 \pmod{p_{d,j}},
    & m<j\le d.
\end{cases}
\]
By the Chinese remainder theorem,
these congruences uniquely determine
the integers $b_1,\ldots,b_{k_d},t$
modulo $n_d$.
We choose their representatives in
$\{1,\ldots,n_d-1\}$.
By construction, all these integers
are coprime to $n_d$.

We claim that the congruence
\[
b_1^{x_1}\cdots b_{k_d}^{x_{k_d}}
\equiv t\pmod{n_d},
\qquad x\in\bits^{k_d},
\]
has a solution if and only if the
original X3C instance has a solution.

For each $1\le j\le m$, put
\[
c_j=\bigl|\Set{i}{1\le i\le\ell,\
x_i=1,\ s_j\in S_i}\bigr|.
\]
Reducing the subset sum congruence
modulo $p_{d,j}$, we obtain
\[
a_j^{c_j}\equiv a_j\pmod{p_{d,j}}.
\]
Since $a_j$ has multiplicative order
$p_{d,j}-1$, this congruence is equivalent to
\[
c_j\equiv1\pmod{p_{d,j}-1}.
\]
Moreover,
\[
0\le c_j\le\ell\le d<p_{d,j}-1.
\]
Hence, the congruence is equivalent
to the integer equality $c_j=1$.

Thus, the subset sum congruence is
equivalent to the system
\[
c_j=1,\qquad 1\le j\le m,
\]
which states precisely that every
element of $S$ belongs to exactly
one of the selected subsets.
The remaining prime factors of $n_d$
impose no additional restrictions.

It remains to verify that the reduction
can be performed in polynomial time.
By the choice of $d$, the parameter $d$
is polynomially bounded in the size
of the X3C instance.
By Lemma~\ref{le:moduli-computation},
the modulus $n_d$ and its prime factors
can be computed in deterministic polynomial
time in $d$.
The same holds for the computation of
the primitive roots $a_j$ and the integers
$b_1,\ldots,b_{k_d},t$.
Furthermore, each of these integers has
binary length $O(d\ln d)$, and
$k_d=O(d\ln d)$.
Therefore, the reduction is polynomial-time.

This proves $\NP$-hardness and hence
$\NP$-completeness.
Since the prime factorization of $n_d$
is computed as part of the reduction,
the problem remains $\NP$-complete
when this factorization is included
in the input.
\end{proof}

Theorem~\ref{thm:ssp-np-complete}
establishes worst-case hardness of
the decision version of $\SSP$.
It does not imply hardness of the
planted subset sum problem
$\SSP_{G_d,k_d}$, whose instances
are always solvable.
The average-case hardness of the
latter problem is established independently
through a reduction from Group-SIS.

Despite the $\NP$-completeness of $\SSP$
over the family $(G_d)$, its positive
instances with $k_d$ increments form
a strongly generic set with respect to
the stratification by the parameter $d$
(see~\cite{GMMU2007}).
Indeed, by Corollary~\ref{cor:solvability}
and Lemma~\ref{lem:units-distance},
the proportion of negative instances
in $G_d^{k_d+1}$ is at most
$\exp(-\Omega(d\ln d))$.
Hence, the constant-time algorithm
that always answers YES is correct
on a strongly generic set of inputs.
This contrasts sharply with the conditional
average-case hardness of the corresponding
planted search problem established above.

\end{document}